\documentstyle{article}
\newtheorem{Definition}{Definition }[section]
\newtheorem{Theorem}[Definition]{Theorem }
\newtheorem{Lemma}[Definition]{Lemma }
\newtheorem{Proposition}[Definition]{Proposition}

\newcommand{\taut}{{\xi}_E}
\newcommand{\Proj}{{\bf P}}

\newcommand{\Segno}{\addtocounter{subsection}{1}
{\noindent \par \bf (\thesubsection)}}
\def\normalbaselines{\baselineskip12pt \lineskip3pt\lineskiplimit3pt\noindent}
\newcommand{\mapright}[1]{\smash{\mathop{\longrightarrow}\limits^{#1}}}
\newcommand{\mapdown}[1]{\Big \downarrow
\rlap{$\vcenter{\hbox{$\scriptstyle#1$}}$}}
\def\proof{\noindent{\bf Proof. }}

\def\claim{\noindent{\bf Claim }}
\def\endproof{\par\hfill $\Box$\par}
\def\O{{\cal O}}

\def\F{{\cal F}}

\def\P{{\bf P}}
\def\Q{{\bf Q}}
\def\P{{\bf P}}

\def\Z{{\bf Z}}

\def\R{{\bf R}}
\def\f{\varphi}
\def\ra{\rightarrow}
\def\iso{\simeq}

\title{ Contractions on a manifold polarized by an ample vector bundle}
\author{M. Andreatta - M. Mella\\Dipartimento di
Matematica,Universit\'a di Trento,\\38050 Povo (TN), Italia\\
e-mail :	Andreatta or Mella @itnvax.science.unitn.it
}
\begin{document}
\normalbaselines
\maketitle

\section*{Introduction}

An algebraic variety $X$ of dimension $n$
(over the complex field) together with an ample vector bundle $E$
on it will be called a {\sf generalized polarized variety}.
The adjoint bundle of the pair $(X,E)$
is the line bundle $K_X + det(E)$.
Problems concerning adjoint bundles have drawn a lot of attention
to algebraic geometer: the classical case is when  $E$
is a (direct sum of) line bundle (polarized variety), while the generalized
case was
motivated by the solution of Hartshorne-Frankel conjecture by Mori (
\cite{Mo})
and by consequent conjectures of Mukai (\cite{Mu}).
\par
A first point of view is to study the positivity (the nefness or
ampleness) of the adjoint line bundle in the case $r = rank (E)$
is about $n = dim X$.
This was done in a sequel of papers for $r\geq (n-1)$
and for smooth manifold $X$ ([Ye-Zhang], [Fujita],
[Andreatta-Ballico-Wisniewski]).
In this paper we want to discuss the next case, namely when $rank (E) = (n-2)$,
with $X$ smooth; we obtain a complete answer which is described in the
theorem (4.1). This is divided in three cases, namely when
$K_X + det(E)$ is not nef, when it is nef and not big and
finally when it is nef and big but not ample.
If $n=3$ a complete picture is already contained in the famous paper
of Mori (\cite{Mo1}), while the particular case in which
$E = \oplus^{(n-2)} (L)$
with $L$ a line bundle was also studied (\cite{Fu1}, \cite{So};
in the singular case see \cite{An}). The part 1 of the theorem was proved
(in a slightly weaker form) by Zhang (\cite{Zh}) and, in the case $E$ is
spanned
by global sections, by Wisniewski (\cite{Wi2}).
\par
Another point of view can be the following: let $(X,E)$ be
a generalized polarized variety with $X$ smooth and $rankE=r$.
If $K_X + det(E)$ is nef, then by the Kawamata-Shokurov
base point free theorem it supports a contraction (see (1.2));
i.e. there exists a map $\pi :X \ra W$ from $X$ onto a
normal projective variety $W$ with connected fiber and such that
$K_X + det(E) = \pi^*H$ for some ample line bundle $H$ on $W$.
It is not difficult to see that, for every fiber $F$ of
$\pi$, we have $dimF \geq (r-1)$, equality holds only if
$dimX > dimW$. In the paper we study the "border" cases:
we assume that $dimF = (r-1)$ for every fibers
and we prove that $X$ has a $\P^r$-bundle
structure given by $\pi$ (theorem (3.2)). We consider also the case
in which $dimF = r$ for every fibers and $\pi$ is birational,
proving that $W$ is smooth and that $\pi$ is a blow-up of a
smooth subvariety (theorem (3.1)).
This point of view was discussed in the case $E = \oplus^r L$
in the paper [A-W].
\par
Finally in the section (4) we extend the theorem (3.2) to the
singular case, namely for projective variety $X$
with log-terminal singularities.  In particular this gives
the Mukai's conjecture1 for singular varieties.

\bigskip

\section{Notations and generalities}
\Segno
 We use the standard notations from algebraic geometry.
Our language is compatible with that of [K-M-M] to which we refer constantly.
We just explain some special definitions
and propositions used frequently.

In particular in this paper $X$ will always stand
for a smooth complex projective variety
of dimension $n$. Let
$Div(X)$ the group of Cartier divisors on $X$;
denote by $K_X$ the {\sf canonical divisor} of $X$, an element of
$Div(X)$ such that $\O_{X}(K_X) = \Omega^n_{X}$.
Let $N_1(X)=\frac{\{1-cycles\}}{\equiv}\otimes \R$,
$N^1(X)= \frac{\{divisors\}}{\equiv}\otimes \R$ and
$\overline {<NE(X)>}=\overline{\{\mbox{effective 1-cycles}\}}$;
 the last is a closed
cone in $N_1(X)$. Let also
$\rho(X)=dim_{\R}N^1(X)<\infty$.

\medskip
Suppose that $K_X$ is not nef, that is there exists an effective curve $C$
such that $K_X\cdot C<0$.

\begin{Theorem}\cite{KMM}
Let $X$ as above and $H$ a nef Cartier divisor such that
$F:= H^{\bot} \cap \overline {<NE(X)>} \setminus \{0\}$
is entirely contained in the set
$\{Z\in N_1(X) :K_X\cdot Z<0\}$,
where $H^{\bot} = \{Z:H\cdot Z=0\}$.
Then there exists a projective morphism $\f:X\ra W$ from $X$ onto a normal
variety $W$ with the following properties:
\begin{itemize}
\item[{i})]  For an irreducible curve $C$ in $X$, $\f(C)$
is a point if and only if $H.C = 0$, if and only if
$cl(C) \in F$.
\item[{ii})] $\f$ has only connected fibers
\item[{iii})] $H = \f^*(A)$ for some ample divisor $A$
on $W$.
\item[{iv})] The image $\f^* :Pic(W) \ra Pic(X)$ coincides with
$\{D \in Pic(X): D.C = 0 \mbox{ \rm  for all } C \in F\}.$
\end{itemize}
\label{contractionth}
\end{Theorem}

\begin{Definition} The following terminology is mostly used (\cite{KMM},
definition 3-2-3).
Referring to the above theorem,

the map $\f$ is called a {\sf contraction}
(or an
{\sf extremal contraction}); the set $F$ is an {\sf extremal face},
while the Cartier
divisor $H$ is a {\sf supporting divisor} for the map $\f$ (or the face $F$).
If $dim_{\R}F = 1$ the face $F$ is called an {\sf extremal ray}, while $\f$
is called an {\sf elementary contraction}.
\end{Definition}

\begin{remark} We have also (\cite{Mo1}) that if $X$ has an extremal ray $R$
then there exists a rational curve $C$ on $X$ such that
$0< -K_X \cdot C\leq n+1$ and
$R=R[C]:=\{D\in <NE(X)>: D\equiv \lambda C, \lambda\in \R^+\}$.
Such a curve is called an {\sf extremal curve}.
\end{remark}

\begin{remark}\label{biraz} Let $\pi:X\ra V$ denote a contraction of an
extremal
face $F$, supported by $H=\pi^*A$([iii]\ref{contractionth}) . Let
 $R$ be an extremal ray in $F$ and $\rho:X\ra W$ the contraction of $R$.
Since $\pi^*A\cdot R=0$, $\pi^*A$ comes from $Pic (W)$
([iv]\ref{contractionth}). Thus $\pi$ factors trough $\rho$.
\end{remark}

\begin{Definition} To an extremal ray $R$ we can associate:
\begin{itemize}
\item[{i})] its {\sf length} $l(R):=min\{ -K_X\cdot C;$
for $C$ rational curve and $C\in R\}$
\item[{ii})] the {\sf locus} $E(R):=\{$the locus
of the curves whose numerical classes are in $R\}\subset X$.
\end{itemize}
\end{Definition}

\begin{Definition} It is usual to divide the elementary contractions
associated to an
extremal ray $R$ in three types according to the dimension of $E(R)$:
more precisely we say that $\f$ is of {\sf fiber type}, respectively
{\sf divisorial type}, resp. {\sf flipping type}, if
$dim E(R) = n$, resp. $n-1$, resp. $< n-1$.
Moreover an extremal ray is said not nef if
there exists an effective  $D\in Div(X)$ such that $D\cdot C<0$.
\end{Definition}

The following very useful inequality was proved in \cite{Io} and \cite{Wi3}.

\begin{Proposition} Let $\varphi$ the contraction of an extremal ray
$R$, $E^{\prime}(R)$ be any irreducible component of the exceptional
locus and $d$ the dimension of a fiber of the contraction
restricted to $E^{\prime}(R)$. Then
$$ dim E^{\prime}(R)+d\geq n+l(R)-1.$$
\label{diswis}
\end{Proposition}

\Segno Actually it is very useful to understand when a contraction is
elementary
or in other words when the locus of two distinct extremal rays are
disjoint. For this we will use in this paper the following results.

\begin{Proposition}\cite[Corollary 0.6.1]{BS} Let $R_1$ and $R_2$ two distinct
not nef
extremal rays such that $l(R_1)+l(R_2)>n$.
Then $E(R_1)$ and $E(R_2)$ are disjoint.
\label{birelementare}
\end{Proposition}

Something can be said also if $l(R_1)+l(R_2)=n$:
\begin{Proposition}\cite[Theorem 2.4]{Fu3}
Let $\pi:X\ra V$ as above and suppose $n\geq 4$ and $l(R_i)\geq n-2$.
Then the exceptional loci corresponding to different extremal rays,
are disjoint with each other.
\label{n=4}
\end{Proposition}

\begin{Proposition}\cite{ABW1} Let $\pi:X\ra W$ be a contraction of
a face such that $dimX > dim W$. Suppose that for every rational curve $C$
in a general fiber of $\pi$ we have
$-K_X\cdot C\geq (n+1)/2$
Then $\pi$ is an elementary contraction except if
\begin{itemize}
\item[a)] $-K_X\cdot C=(n+2)/2$ for some rational curve $C$ on $X$,
$W$ is a point, $X$ is a Fano manifold
of pseudoindex $(n+2)/2$ and $\rho(X)=2$
\item[b)]  $-K_X\cdot C=(n+1)/2$
for some rational curve $C$, and $dim$W$\leq 1$
\end{itemize}
\label{fibelementare}
\end{Proposition}

The following definition is used in the theorem:
\begin{Definition} Let $L$ be an an
ample line bundle on $X$. The pair $(X,L)$ is called a scroll (respectively
a quadric fibration, respectively a del Pezzo fibration)
over a normal variety
$Y$ of dimension $m$ if there exists a surjective morphism
with connected fibers
$\phi: X \ra Y$ such that
$$K_X+(n-m+1)L \approx p^*{\cal L}$$
(respectively $K_X+(n-m)L \approx p^*{\cal L}$;
respectively $K_X+(n-m-1)L \approx p^*{\cal L}$)
for some ample line bundle ${\cal L}$ on $Y$.
$X$ is called a classical scroll (respectively
quadric bundle) over a projective
variety $Y$ of dimension $r$ if there exists a surjective morphism
$\phi : X\ra Y$ such that every fiber is isomorphic to $\P^{n-r}$
(respectively to a quadric in $\P^{(n-r+1)}$)
and if there exists a vector bundle $E$ of rank $(n-r+1)$ (respectively
of rank $n-r+2$) on $Y$ such that  $X\simeq \P(E)$
(respectively exists an embedding of $X$ as a subvariety of $\P(E)$).
\end{Definition}

\section{A technical construction}
\label{tech}
Let $E$ be a vector bundle of rank $r$ on $X$ and assume
that $E$ is ample, in the sense of Hartshorne.

\begin{remark} Let $f:{\bf P}^1\rightarrow X$
be a non constant map, and $C=f({\bf P}^1)$,
\label{mha}
then $detE\cdot C\geq r$.
\par
In particular if there exists a curve $C$ such that
$(K_X+detE).C \leq 0$ (for instance if $(K_X+detE)$
is not nef) then there exists
an extremal ray $R$ such that  $l(R) \geq r$.
\end{remark}

\Segno
\label{sopra}
Let $Y=\Proj(E)$ be the associated projective space bundle,
$p:Y \rightarrow X$ the natural map onto $X$ and
$\taut$ the tautological bundle of $Y$.
Then we have the formula for the canonical bundle
$K_Y=p^*(K_X+detE)-r\taut$.
Note that $p$ is an elementary contraction; let $R$ be the associated extremal
ray.

Assume that $K_X+detE$ is nef but not ample and that it is
the supporting divisor of an elementary contraction $\pi:X\rightarrow W$.
Then $\rho(Y/W) = 2$ and $-K_Y$ is $\pi\circ p$-ample.
By the relative Mori theory
over $W$ we have that there exists a ray on $NE(Y/W)$, say
$R_1$, of length $\geq r$, not contracted by $p$, and a relative
elementary contraction
$\f:Y\ra V$. We have thus the following commutative diagram.

\begin{equation}
\label{dia1}
\matrix{{\bf P}(E)=Y&\mapright\varphi&V\cr
\mapdown{p}&&\mapdown\psi\cr X&\mapright\pi&W}
\end{equation}
where $\varphi$ and $\psi$ are elementary contractions.
Let $w\in W$ and let $F(\pi)_w$ be an irreducible component of $\pi^{-1}(w)$;
choose
also $v$ in $\psi^{-1}(w)$ and let $F(\varphi)_v$ be an irreducible
component of $\varphi ^{-1}(v)$
such that $p(F(\varphi)_v) \cap F(\pi)_w \not= \emptyset$; then
$p(F(\varphi)_v) \subset F(\pi)_w$. This is true by the commutativity of the
diagram.
Since $p$ and $\f$ are elementary contractions of different extremal
rays we have that $dim(F(\varphi)\cap F(p))=0$,
that is curve contracted by $\varphi$ cannot be contracted by $p$.

In particular this implies that $dim p(F(\varphi)_v) = dim F(\varphi)_v$;
therefore
$$dimF(\varphi)_v\leq dimF(\pi)_w.$$

\begin{remark} If $dimF(\varphi)_v=dimF(\pi)_w$,
then $dim F(\psi)_w:=dim(\psi^{-1}(w)) = r-1$;
if this holds for every $w \in W$ then $\psi$ is equidimensional.
\end{remark}

\proof Let $Y_w$ be an irreducible component of
$p^{-1}\pi^{-1}(w)$ such that $\varphi (Y_w) = F(\psi)_w$.
Then $dim F(\psi)_w = dim Y_w - dimF(\varphi)_v = dim Y_w -
dimF(\pi)_w = dimF(p) = (r-1)$.
\endproof

\Segno {\bf Slicing techniques}
\label{adj}

Let $H = \varphi ^*(A)$ be a supporting divisor for $\varphi$ such that
the linear system $|H|$ is base point free.
We assume as in (\ref{sopra}) that $( K_X + detE)$ is nef
and we refer to the diagram (\ref{dia1}). The divisor
$K_Y+r\taut =p^*( K_X + detE)$ is nef on $Y$ and therefore
$m(K_Y+r\taut+aH)$, for $m\gg 0$, $a\in{\bf N}$,
is also a good supporting divisor for $\varphi$.
Let $Z$ be a smooth n-fold obtained by intersecting $r-1$ general divisor
from the linear system H, i.e. $Z = H_1\cap \dots \cap H_{r-1}$ (this is what
we call a {\sf slicing});
let $H_i = \varphi^{-1} A_i$.

Note that the map ${\varphi}^{\prime}=\varphi_{|Z}$
is supported by $m|(K_Y+r\taut+a\varphi^*A)_{|Z}|$,
hence, by adjunction, it is supported by $K_Z+rL$, where $L={\taut}_{|Z}$.
Let $p^{\prime}=p_{|Z}$; by construction $p^{\prime}$ is finite.

If $T$ is (the normalization of) $\varphi (Z)$ and
$\psi^{\prime} :T \ra W$ is the map obtained restricting $\psi$
then we have from (\ref{dia1}) the following diagram
\begin{equation}
\label{dia2}
\matrix{Z&\mapright{\varphi\prime}&T\cr
\mapdown{p\prime}&&\mapdown{\psi\prime}\cr X&\mapright\pi&W}
\end{equation}

In general one has a good comprehension of the map $\f^{\prime}$
(for instance in the case $r = (n-2)$ see the results in \cite{Fu1} or in
\cite{An}). The goal is to "transfer" the information that we have on
$\f^{\prime}$ to the map $\pi$. The following proposition
is the major step in this program.

\begin{Proposition} Assume that
$\psi$ is equidimensional (in particular this is the case if
for every non trivial fiber we have
$dimF(\varphi)=dimF(\pi)$). Then $W$ has the same singularities of $T$.
\label{fujita}
\end{Proposition}

\begin{proof} By hypothesis any irreducible reduced component $F_i$ of a  non
trivial
fiber $F(\psi)$
is of dimension $r-1$; this implies also that
$F_i=\varphi(F(p))$ for some fiber of $p$.

Now, let us follow an argument as in
\cite[Lemma 2.12]{Fu1}.
We can assume that the divisor $A$ is very ample; we will choose
$r-1$ divisors $A_i \in |A|$ as above such that, if
$T = {\bigcap_{i}} A_i$, then $T \cap \psi^{-1}(w)_{red}= N$
is a reduced 0-cycle and $Z = H_1\cap \dots \cap H_{r-1}$ is a
smooth n-fold, where $H_i = \varphi^{-1} A_i$. This can be done by
Bertini theorem. Moreover the number of points in $N$ is given by
$A^{r-1}\cdot \psi^{-1}(w)_{red}=\sum_i A^{r-1}\cdot F_i=\sum_i d_i$.
Note that, by projection formula, we have
$A^{r-1}\cdot F_i= \varphi^*A^{r-1}\cdot F(p)$;
moreover, since $p$ is a projective bundle, the last number is constant
i.e. $\varphi^*A^{r-1}\cdot F(p) = d$ for all fiber $F(p)$,
that is the $d_i$'s  are constant.

Now take a small enough neighborhood $U$ of $w$, in the metric topology,
such that any
connected component $U_{\lambda}$ of $\psi^{-1}(U)\cap T$
meets $\psi^{-1}(w)$ in a single point. This is
possible because $\psi^{\prime}:=\psi_{|T}: T \ra W$
is proper and finite over $w$.
Let $\psi_{\lambda}$ the restriction of $\psi$ at
$U_{\lambda}$ and $m_{\lambda}$ its degree.
Then $deg{\psi}^{\prime}=\sum m_{\lambda}\geq \sum_i d_i = \sum_i d$
and equality holds if and only if $\psi$ is not ramified at $w$
(remember that $\sum_i d_i$ is the number of $U_{\lambda}$).

The generic $F(\psi)_w$ is irreducible and generically
reduced.
Note that we can choose
$\tilde{w}\in W$ such that $\psi^{-1}(\tilde{w})=
\varphi(F(p))$ and
$deg{\psi}^{\prime}=A^{r-1}\cdot\psi^{-1}(\tilde{w})$, the latter is
possible by  the choice of generic sections of $|A|$.
Hence, by projection formula
$deg\psi^{\prime}= A^{r-1}\cdot \psi^{-1}(\tilde{w})=
\varphi^*A^{r-1}\cdot F(p)=d$,
that is $m_{\lambda}=1$ and the fibers are irreducible.
Since $W$ is normal we can conclude, by Zarisky's Main theorem, that $W$ has
the
same singularity as $T$.
\endproof

\section{Some general applications}

As an application of the above construction we will prove the following
proposition; the case $r = (n-1)$ was proved in
\cite {ABW2}.

\begin{Proposition} Let $X$ be a smooth projective complex variety and $E$ be
an ample vector bundle of rank $r$ on $X$. Assume that $K_X+detE$ is nef and
big but
not ample and let $\pi:X\ra W$ be the
contraction supported by $K_X+detE$.
Assume also that $\pi$ is a divisorial elementary contraction,
with exceptional divisor $D$,
and that $dim F\leq r$ for all fibers $F$.
Then $W$ is smooth, $\pi$ is the blow up of a
smooth subvariety $B: = \pi (D)$ and
$E =\pi^*E^{\prime}\otimes[-D]$, for
some ample $E^{\prime}$ on $W$. \label{bd}
\end{Proposition}

\proof Let $R$ be the extremal ray contracted by $\pi$ and $F:=F(\pi)$ a fiber.
Then $l(R)\geq r$ and thus $dimF\geq r$ by proposition (\ref{diswis}).
Hence all the fibers of $\pi$ have dimension $r$.
Consider the commutative diagram (\ref{dia1}); let $R_1$ be
the ray contracted by $\f$. Since $l(R_1)\geq r$, again
by proposition (\ref{diswis}), we have that $dimF(\f) \geq r$ (note that
$R_1$ is not nef). Therefore,
since $dimF(\f) \leq  dimF$, we have that
$dimF(\f) =  dimF = r$, $l(R)=l(R_1)=r$ and $\taut\cdot C_1=1$, where $C_1$ is
a (minimal) curve in the ray $R_1$.
Via slicing we obtain the map $\varphi^{\prime}:Z\ra T$ which is supported by
$K_Z+r\taut{}_{|Z}$. This last map is very well understood:
namely by \cite[Th 4.1 (iii)]{AW} it follows that $T$ is smooth and
$\varphi^{\prime}$ is a blow up along a smooth subvariety.
By proposition (\ref{fujita}) also $W$ is smooth.
Therefore $\pi$ is a birational morphism between smooth varieties
with exceptional locus a prime divisor and with equidimensional
non trivial fibers; by \cite[Corollary 4.11]{AW} this implies that
$\pi$ is a blow up of a smooth subvariety in $W$.

We want to show  that $E =\pi^*E^{\prime}\otimes[-D]$.
Let $D_1$ be the exceptional divisor of $\varphi$; first we claim that
$\taut+D_1$ is a good supporting divisor for $\f$.
To see this observe that $(\taut+D_1)\cdot C_1=0$, while $(\taut+D_1)\cdot C>0$
for any curve $C$ with $\f(C)\not= pt$ (in fact $\taut$ is ample and $D_1\cdot
C\geq 0$ for such a curve).
Thus $\taut+D_1=\f^*A$ for some ample $A\in Pic(V)$; moreover by
projection formula $A\cdot l=1$, for any line $l$ in the fiber of $\psi$.
Hence by Grauert theorem $V=\Proj(E^{\prime})$ for some ample vector bundle
$E^{\prime}$ on $W$. This yields, by the commutativity of diagram (1), to
$E\otimes D=p_*(\taut+D_1)=p_*\f^*A=\pi^*\psi_*A=\pi^*E^{\prime}$.
\endproof

We now want to give a similar proposition for the fiber type case.

\begin{Theorem} Let $X$ be a smooth projective complex variety and $E$
be an ample vector bundle of rank $r$ on $X$. Assume that $K_X+detE$
is nef and let $\pi:X\ra W$ be the
contraction supported by $K_X+detE$.
Assume that $r\geq (n+1)/2$ and
$dim F\leq r-1$ for any fiber $F$ of $\pi$. Then $W$ is smooth,
for any fiber $F\simeq \Proj^{r-1}$ and $E_{|F}=\oplus^r\O(1)$.
\label{relscroll}
\end{Theorem}

\proof Note that by proposition (\ref{diswis}) $\pi$ is a contraction of
fiber type and all the fibers have dimension $r-1$.
Moreover the contraction is elementar, as it follows from
proposition (\ref{fibelementare}).

We want to use an inductive argument to prove the thesis.
If $dim W=0$ then this is Mukai's conjecture1; it was proved
by Peternell, Koll\'ar, Ye-Zhang (see for instance \cite{YZ}).
Let the claim be true for dimension $m-1$.
Note that the locus over which the fiber is not $\Proj^{r-1}$
is discrete and $W$ has isolated singularities.
In fact take a general hyperplane section $A$ of $W$,
and $X^{\prime}=\pi^{-1}(A)$ then $\pi_{|X^{\prime}}:X^{\prime}\ra A$
is again a contraction supported by $K_{X^{\prime}}+det E_{|X^{\prime}}$,
such that $r\geq ((n-1)+1)/2$. Thus by induction $A$ is smooth, hence $W$
has isolated singularities.
\par
Let $U$ be an open disk in the complex topology, such that $U\cap SingW=\{0\}$.
Then by lemma below \ref{scroll} we have locally, in the complex topology, a
$\pi$-ample line bundle $L$ such that restricted to the general fiber is
$\O(1)$. As in \cite[Prop. 2.12]{Fu1} we can prove
that $U$ is smooth and all the fibers are $\Proj^{r-1}$.
\endproof

\begin{Lemma}
\label{scroll}
Let $X$ be a complex manifold and $(W,0)$ an analityc germ such that
$W\setminus \{0\}\simeq \Delta^m\setminus \{0\}$.
Assume we have an holomorphic map $\pi:X\ra W$ with $-K_X$ $\pi$-ample;
assume also that $F\simeq \Proj^r$ for all fibers of $\pi$, $F\not=
F_0=\pi^{-1}(0)$,
and that $codim F_0\geq 2$.
Then there exists a line bundle $L$ on $X$ such that $L$ is $\pi$-ample and
$L_{|F}=\O(1)$.
\end{Lemma}

\proof (see also \cite[pag 338, 339]{ABW2})
Let $W^*=W\setminus \{0\}$ and $X^*=X\setminus F_0$.
By abuse of notation call $\pi=\pi_{|X^*}:X^*\ra W^*$;
it follows immediately that $R^1\pi_*\Z_{X^*}=0$ and $R^2\pi_*\Z_{X^*}=\Z$.

If we look at Leray spectral sequence, we have that:
$$
E^{0,2}_2= \Z\mbox{ and } E^{p,1}_2= 0 \mbox{ for any p.}$$
Therefore $d_2:E^{0,2}_2\ra E^{2,1}_2$ is the zero map
and moreover we have the following exact sequence
$$0\ra E^{0,2}_{\infty}\ra E^{0,2}_2\stackrel{d_3}{\ra} E^{3,0}_2,$$
since	the only non zero map from $ E^{0,2}_2$ is $d_3$ and hence
$E^{0,2}_{\infty}=kerd_3$.
On the other hand we have also, in a natural way, a surjective map
$H^2(X^*,\Z)\ra  E^{0,2}_{\infty}\ra 0$.
Thus we get the following exact sequence
$$ H^2(X^*,\Z)\stackrel{\alpha}{\ra} E^{0,2}_2\ra
E^{3,0}_2=H^3(W^*,\Z).$$ We want to show that $\alpha$ is surjective.
If $dimW := w\geq 3$ then $H^3(W^*,\Z)=0$ and we have done. Suppose $w=2$ then
$H^3(W^*,\Z)=\Z$; note
that the restriction of $-K_X$ gives a non zero class (in fact it is
$r+1$ times the generator) in $E^{0,2}_2$ and is mapped to zero in $E^{0,3}_2$
thus the mapping
$E^{0,2}_2\ra E^{3,0}_2$ is the zero map and $\alpha$ is surjective.
Since $F_0$ is of
codimension at least 2 in $X$ the restriction map $H^2(X,\Z)\ra
H^2(X^*,\Z)$ is a bijection. By the vanishing of $R_i\pi_*\O_X$ we get
$H^2(X,\O_X)=H^2(W,\O_W)=0$ hence also $Pic(X)\ra H^2(X,\Z)$ is surjective.
Let $L\in Pic(X)$ be a preimage of a generator of $E^{0,2}_2$.
 By construction $L_t$	is $\O(1)$, for $t\in W^*$. Moreover
$(r+1)L=-K_X$ on $X^*$ thus, again by the codimension of $X^*$, this is true
on $X$ and $L$ is $\pi$-ample.
\endproof

\section{An approach to the singular case}
 The following theorem arose during a discussion between us
and J.A. Wisniewski; we would like to thank him.
The idea to investigate this argument came from a preprint of Zhang
[Zh2] where he proves the following result under
the assumption that $E$ is spanned by global sections.
For the definition of log-terminal
singularity we refer to \cite{KMM}.

\begin{Theorem} Let $X$ be an n-dimensional log-terminal projective
variety and $E$ an ample vector bundle of rank $n+1$, such that
$c_1(E)=c_1(X)$.
Then $(X,E)=(\Proj^n,\oplus^{n+1}\O_{\Proj^n}(1))$.
\end{Theorem}

\proof We will prove that $X$ is smooth, then we can apply proposition
(\ref{relscroll}).
We consider also in this case the associated projective space bundle $Y$
and the commutative diagram
\begin{equation}
\matrix{{\bf P}(E)=Y&\mapright\varphi&V\cr
\mapdown{p}&&\mapdown\psi\cr X&\mapright\pi&pt}
\end{equation}
as in (\ref{dia1}); it is immediate that  $Y$ is a weak Fano variety
(i.e. $Y$ is Gorenstein, log-terminal and $-K_Y$ is ample; in particular
it has Cohen-Macaulay singularities); moreover, as in (3.1),
$dimF(\f) \leq dimF(\pi) = n$ and the map $\f$ is supported by $K_Y+(n+1)H$,
 where
$H =\taut + A$, with $\taut$ the tautological line bundle and $A$ a
pull back of a ample line bundle from $V$.
It is known that a contraction supported
by $K_Y+rH$ on a log terminal variety has to have fibers
of dimension $\geq (r-1)$ and of dimension
$\geq r$ in the birational case (\cite [remark 3.1.2]{AW}).
 Therefore in our case
$\f$ can not be birational and all fibers have dimension $n$;
moreover, by the Kobayashi-Ochiai criterion the general fiber is
 $F\iso \Proj^n$.
We want to adapt the proof of \cite[Prop 1.4]{BS};
to this end we have only to show
that there are no fibers of $\f$ entirely contained in $Sing(Y)$. Note
that, by construction, $Sing(Y)\subset p^{-1}(Sing X)$
hence no fibers $F$ of $\f$
can be contained in $Sing(Y)$.
Hence the same proof of \cite[Prop 1.4]{BS} applies and we can prove that
$V$ is nonsingular and $\f:Y\ra V$ is a classical scroll.
In particular $Y$ is nonsingular and therefore also $X$ is nonsingular.
\endproof

More generally we can prove the following.
\begin{Theorem}  Let $X$ be an n-dimensional log-terminal projective
variety and $E$ be an ample vector bundle of rank $r$. Assume that $K_X+det E$
is nef and let $\pi:X\ra W$ be the contraction supported by $K_X+det E$.
Assume also that for any fiber $F$ of $\pi$ $dimF\leq r-1$, that $r\geq
(n+1)/2$
and $codim Sing(X)>dim W$.
Then $X$ is smooth and for any fiber $F\simeq \Proj^{r-1}$.
\end{Theorem}
\proof The proof that $X$ is smooth is as in the theorem above and then
we use proposition (\ref{relscroll})
\endproof

\section{Main theorem}

This section is devoted to the proof of the following theorem.
\begin{Theorem} Let $X$ be a smooth projective variety over the complex
field of dimension $n \geq 3$ and $E$ an ample vector bundle on $X$
of rank $r= (n-2)$. Then we have
\begin{itemize}

\item[1)] $K_X + det(E)$ is nef unless $(X,E)$ is one of the following:
\begin{itemize}
\item[{i})] there exist a smooth $n$-fold, $W$,  and
a morphism $\phi : X \ra W$ expressing $X$ as a blow up of a
finite set $B$ of points and an ample vector bundle $E'$ on
$W$ such that $E = \phi^*E'\otimes[-\phi^{-1}(B)]$.

\par\noindent
Assume from now on that $(X,E)$ is not as in (i) above (that is eventually
consider the new pair $(W,E')$ coming from (i)).
\item[{ii})] $X = \P^n$ and $E =\oplus^{(n-2)}\O(1)$ or
$\oplus^{2}\O(2)\oplus^{(n-4)}\O(1)$ or
$\O(2)\oplus^{(n-3)}\O(1)$ or
$\O(3)\oplus^{(n-3)}\O(1)$.
\item[{iii})] $X = \Q^n$ and $E =\oplus^{(n-2)}\O(1)$ or
$\O(2)\oplus^{(n-3)}\O(1)$ or ${\bf E}(2)$ with ${\bf E}$
a spinor bundle on $\Q^n$.
\item[{iv})] $X = \P^2 \times \P^2$ and $E = \oplus^2\O(1,1)$
\item[{v})] $X$ is a del Pezzo manifold with $b_2 = 1$, i.e.
$Pic(X)$ is generated by an ample line bundle $\O(1)$ such that
$\O(n-1) = \O(-K_X)$ and $E = \oplus^{(n-1)}\O(1)$.
\item[{vi})] $X$ is a classical scroll or a quadric bundle over
a smooth curve $Y$.
\par
\item[{vii})] $X$ is a fibration over a smooth surface
$Y$ with all fibers isomorphic to $\P^{(n-2)}$.
\end{itemize}

\item[2)] If $K_X + det(E)$ is nef then it is big unless
there exists a morphism $\phi : X \ra W$ onto a
normal variety $W$ supported by (a large multiple of) $K_X + det(E)$
and $dim(W) \leq 3$; let $F$ be a general fiber of $\phi$ and
$E^{\prime}=E_{|F}$.
We have the following according to
$s = dim W$:
\begin{itemize}
\item[{i})] If $s = 0$ then $X$ is a Fano manifold and
$K_X + det(E) = 0$. If $n\geq 6$ then $b_2(X) = 1$ except if
$X=\P^3\times\P^3$ and $E=\oplus^4\O(1,1)$.
\item[{ii})] If $s = 1$ then $W$ is a smooth curve and $\phi$ is a
flat (equidimensional) map.
Then $(F,E')$ is one of the pair described in \cite{PSW},
in particular $F$ is either $\P^n$ or a quadric or a del Pezzo variety.
If $n \geq 6$ then $\pi$ is an elementary contraction.
If the general fiber is $P^{n-1}$ then $X$ is a classical scroll while
if the general fiber is $\Q^{n-1}$ then $X$ is a quadric bundle.
\item[{iii})] If $s = 2$ and $n \geq 5$ then $W$ is a smooth surface,
$\phi$ is a flat map and $(F,E^{\prime})$ is one of the pair described
in the Main Theorem of \cite{Fu2}.
If the general fiber is $\Proj^{n-2}$ all the fibers are $\Proj^{n-2}$.
\item[{iv})] If $s = 3$ and $n \geq 5$ then $W$ is a smooth
3-fold and all fibers are isomorphic to $\Proj^{n-3}$.
\end{itemize}
\end{itemize}
\item[3)] Assume finally that $K_X + det(E)$ is nef and big but not ample.
Then a high multiple of $K_X + det(E)$ defines a birational map,
$\f :X \ra X'$, which contracts an "extremal face" (see section 2).
Let $R_i$, for $i$ in a finite set of index, the extremal rays
spanning this face; call
$\rho_i: X \ra W$ the contraction associated to one of the $R_i$.
Then we have that each $\rho_i$ is birational and divisorial; if $D$ is
one of the exceptional divisors
 (we drop the index) and $Z = \rho (D)$
we have that $dim(Z) \leq 1$
and the following possibilities occur:
\begin{itemize}
\item[{i})] $dimZ = 0$, $D = \P^{(n-1)}$ and $D_{|D} = \O(-2)$
or $\O(-1)$;
moreover, respectively, $E_{|D} =\oplus^{n-2}\O(1)$ or
 $E_{|D} =\oplus^{n-1}\O(1)\oplus \O(2)$.
\item[{ii})] $dimZ = 0$, $D$ is a (possible singular) quadric,
$\Q^{(n-1)}$, and $D_{|D} = \O(-1)$;
moreover  $E_{|D} =\oplus^{n-2}\O(1)$.
\item[{iii})] $dimZ = 1$, $W$ and $Z$ are smooth projective varieties
and $\rho$ is the blow-up of $W$ along $Z$.
Moreover $E_{|F} =\oplus^{n-2}\O(1)$.
\end{itemize}
If $n > 3$ then $\f$ is a composition of "disjoint" extremal
contractions as in i), ii) or iii).
\label{main}
\end{Theorem}

\proof Proof of part 1) of the theorem

Let $(X,E)$ be a generalized polarized variety
and assume that $K_X + det(E)$ is not nef.
Then there exist on $X$ a finite number of extremal rays, $R_1, \dots , R_s$,
such that $(K_X + det(E))^.R_i < 0$ and therefore, by the remark in section
(2),
$l(R_i) \geq (n-1)$.

Consider one of this extremal rays, $R = R_i$, and let $\rho : X \ra Y$
be its associated elementary contraction. Then $L := -(K_X+det(E))$ is
$\rho$-ample and also the vector bundle $E_1 := E \oplus L$ is
$\rho$-ample; moreover $K_X + det(E_1) = \O_X$ relative to $\rho$.
We can apply the theorem in \cite{ABW2} which study the positivity of
the adjoint bundle in the case of $rank E_1 = (n-1)$. More precisely
we need a relative version of this theorem, i.e. we do not assume
that $E_1$ is ample but that it is $\rho$-ample
(or equivalently a local statement in a neighborhood
of the exceptional locus of the extremal ray $R$).
We just notice that the theorem in \cite{ABW2} is true also in the relative
case
and can be proved exactly with the
same proof using the relative minimal model theory (see [K-M-M]; see also the
section 2
of the paper \cite{AW} for a discussion of the local set up).

Assume first that $\rho$ is birational, then $K_X + det(E_1)$ is $\rho$-nef and
$\rho$-big; note also that, since $l(R_i) \geq (n-1)$, $\rho$ is divisorial.
Therefore we are in the (relative) case C of the theorem
in \cite{ABW2} (see also the proposition \ref{bd}
with $r = (n-1)$); this implies that $Y$ is smooth and $\rho$ is the blow up
of a point in $Y$.
Since $l(R_i) \geq (n-1)$, the exceptional loci
of the birational rays are pairwise disjoint by proposition
(\ref{birelementare}).
This part give the point {\sf (i)} of the
theorem \ref{main}; i.e. the birational extremal rays have
disjoint exceptional loci which are divisors isomorphic
to $\P^{(n-1)}$ and which contract
simultaneously to smooth distinct points on a $n$-fold $W$. The
description of $E$ follows trivially (see also \cite{ABW2}).

If $\rho$ is not birational then we are in the case
B of the theorem in \cite{ABW2}; from this we
obtain similarly as above the other cases of the theorem \ref{main},
 with some trivial
computations needed to recover $E$ from $E_1$.

\endproof

\bigskip
Proof of the part 2) of the theorem

Let $K_X+detE$ be nef but not big; then it is the supporting
divisor of a face $F = (K_X+detE)^{\bot}$. In particular we can
apply the theorems of section (\ref{tech}): therefore there
exist a map $\pi:X\rightarrow W$ which is given by a high multiple
of $K_X+detE$ and which contracts the curves in the face. Since
$K_X+detE$ is not big we have that  $dimW<dimX$.
 Moreover for every rational curve $C$
in a general fiber of $\pi$ we have $-K_X\cdot C \geq (n-2)$ by the remark
in section (\ref{tech}).
We apply proposition
(\ref{fibelementare}), which, together with the above inequality on
$-K_X\cdot C$, says that $\pi$ is an elementary contraction
if $n\geq 5$ unless
either $n=6$, $W$ is a point and $X$ is a Fano manifold of pseudoindex
$4$ and $\rho(X) = 2$ or $n = 5$ and $dimW \leq 1$.

By proposition (\ref{diswis}) we have the inequality
$$n+dimF\geq n+n-2-1;$$
in particular it follows that $dim W\leq 3$.

\Segno Let $dimW=0$, that is $K_X+detE=0$ and therefore $X$ is a Fano manifold.
By what just said above we have that $b_2(X)=1$ if $n \geq 6$ with an exception
which will be treated in the following lemma.

\begin{Lemma} Let $X$ be a $6$ dimensional projective manifold,
$E$ is an ample vector bundle on $X$ of rank $4$ such that $K_X+detE=0$.
Assume moreover that $b_2 \geq 2$.
Then $X=\P^3\times\P^3$ and $E=\oplus^4\O(1,1)$.
\label{slice}
\end{Lemma}

\proof
The lemma is a slight generalizzation of \cite[Prop B]{Wi1} for
dimension $6$; the poof is similar
and we refer to this paper. In particular as in \cite {Wi1}
we can see that $X$ has two extremal rays whose contractions,  $\pi_i$,$i
=1,2$,
are of fiber type with equidimensional fibers onto
3-folds $W_i$ and with general fiber $F_i\simeq \P^3$.
We claim that the $W_i$ are smooth and thus $W_i\simeq \P^3$.
First of all note that $W_i$ can have only isolated singularity and only
isolated points over which the fiber is not $\Proj^{n-3}$; in fact
let $S$ be a general hyperplane section of $W_i$ and $T_i=\pi_i^*(S)$, then
$(\pi_i)_{|T_i}$ is an extremal contraction, by proposition
\ref{fibelementare};
hence by \cite[Prop 1.4.1]{ABW2} $S$ is smooth; moreover the contraction is
supported by $K_{T_i}+det E_{T_i}$ hence all fibers are
$\P^3$ by the main theorem of \cite{ABW2}.
Now we are (locally) in the hypothesis of lemma \ref{scroll} so we get,
locally in the complex topology, a tautological bundle and
we can conclude, by \cite[Prop 2.12]{Fu1}, that $W_i$ is smooth.
Let $T = H_1 \cap H_2$, where $H_i$ are two general
elements of $\pi_1^*(\O(1)$. $T$ is smooth, we claim that
$T\simeq \Proj^1\times \Proj^3$. In fact $\pi_{1 _{|T}}$
makes $T$ a projective bundle over a line (since $H^2(\Proj^1,\O^*)=0$),
that is $T=\Proj(\F)$. Moreover $\pi_{2_{|T}}$ is onto $\Proj^3$,
therefore the claim follows. Therefore we conclude
that $\pi_i^*\O_{\Proj^3}(1)_{|F_i}\simeq \O_{\Proj^3}(1)$ for $i=1,2$.
This implies by Grauert Theorem that the two fibrations are classical scroll,
 that is
$X=\P(\F_i)$, for $i=1,2$; moreover
computing the canonical class of $X$ the $\F_i$ are ample
and the lemma easily follows.
\end{proof}

\Segno Let $dimW=1$. Then $W$ is a smooth curve and $\pi$ is a flat map.
Let $F$ be a general
fiber, then $F$ is a smooth Fano manifold and $E_{|F}$ is
an ample vector bundle on $F$ of rank $(n-2) = dimF - 1$
such that $-K_F = det(E_{|F})$. These pairs $(F, E_{|F})$
are classified in the Main Theorem of \cite{PSW}; in particular if $dimF \geq
5$
$F$ is either $\P^{(n-1)}$ or $\Q^{(n-1)}$ or a
del Pezzo manifold with $b_2(F) = 1$.
Moreover if $n \geq 6$ then $\pi$
is an elementary contraction by proposition (\ref{fibelementare}).

\claim Let $n\geq 6$ and assume that the general fiber is $\P^{n-1}$, then
 $X$ is a classical scroll and $E_{|F}$ is the same for all $F$.

(See also \cite {Fu2}) Let $S= W\setminus U$ be the locus
of points over which the fiber is not $\Proj^{n-1}$. Over $U$ we have  a
projective
fiber bundle. Since $H^2(U,{\cal O}^*)=0$ we
can associate this $\P$-bundle to a vector bundle $\F$ over $U$.
Let $Y=\P(\F)$ and
$H$ the tautological bundle;
by abuse of language let $H$ the extension of $H$ to $X$. Since
$\pi$ is elementary $H$ is an ample line bundle on $X$.
Therefore by semicontinuity $\Delta(F,H_F)\geq \Delta(G,H_G)$, for
any fiber $G$, where $\Delta(X,L)$ is Fujita delta-genus. In our case
this yields $0=\Delta(F,H_F)\geq \Delta(G,H_G)\geq 0$. Moreover by
flatness $(H_G)^{n-1}=(H_F)^{n-1}=1$ and Fujita classification allows to
conclude.
The possible vector bundle restricted to the fibers are all decomposables,
hence they are rigid, that is $H^1(End(E))=\oplus_i H^1(End(\O(a_i))=
\oplus_i H^1(\O(-a_i))=0$. Hence the decomposition is the same
 along all fibers of $\pi$.

\claim Let $n\geq 6$ and assume that the general fiber is $\Q^{n-1}$.
Then $X$ is a quadric bundle.

Let as above $S=W\setminus U$ be the locus of points over which the fiber is
not
a smooth quadric.
Let $X^*=\pi^{-1}(U)$ then we can embed $X^*$ in a fiber bundle of projective
spaces over $U$,
since it is locally trivial. Associate this $P$-bundle over $U$
to a projective bundle and argue as before.
\endproof

\Segno Let now $dimW=2$ and assume that $n\geq 5$;
then $\pi$ is an elementary contraction.
This implies first, by
\cite[Prop. 1.4.1]{ABW2}, that $W$
is smooth; secondly that $\pi$ is equidimensional, hence flat and
the general fiber
is $\Proj^{n-2}$ or ${\bf Q}^{n-2}$, see \cite{Fu2}.

\claim Let $n\geq 5$ and  the general fiber is $\P^{n-2}$ then
 for any fiber $F\simeq \Proj^{n-2}$ and
$E_{|F}$ is the same for all $F$.

Let $S\subset W$
be the locus of singular fibers, then $dimS\leq 0$ since
$W$ is normal. Let $U\subset W$
be an open set, in the complex topology, with $U\cap S=\{0\}$ and let
$V\subset X$ such that $V=\pi^{-1}(U)$. We are in the hypothesis
of lemma \ref{scroll}
thus we get a  "tautological" line bundle $H$ on
$V$ and we conclude by \cite[Prop. 2.12]{Fu1}.

There are two possible restriction of $E$ to the fiber,
namely $E_{|F}\simeq \O(2)\oplus(\oplus^{n-1}\O(1))$ or
$E_{|F}$ is the tangent bundle. As observed by Fujita in \cite{Fu2} this two
restrictions have a different behavior in the diagram (\ref{dia1}),
in the former $\varphi$ is birational while in the latter it is of fiber
type. Hence the restriction has to be constant along all the fibers.
\endproof

\Segno Let finally $dimW=3$; the general fiber is $\Proj^{n-3}$
(see for instance \cite{Fu2}).
 Assume that $n\geq 5$, therefore
$\pi$ is elementary; we claim that all fibers are $\Proj^{n-3}$.

Since $\pi$ is elementary any fiber $G$ has $cod G\geq 2$.
Let $S\subset W$
be the locus of point over which the fiber is not $\Proj^{n-3}$;
$dimS\leq 0$ since a generic linear
space section can not
intersect $S$, by the above.
Let $U\subset W$
be an open set, in the complex topology, with $U\cap S=\{0\}$ and let
$V\subset X$ such that $\pi(V)=U$.
Then by lemma \ref{scroll} we get a  "tautological" line bundle $H$ on
$V$; $\pi: V\ra U$  is supported by $K_V+(n-2)H$.
Thus by  \cite[Th 4.1]{AW}  $U$ is smooth and all the fibers are $\Proj^{n-3}$
( we use that $n\geq 5$).
\endproof

Proof of the part 3) of the theorem

In the last part of the theorem we assume that $K_X+detE$ is nef and big but
not ample.
Then $K_X+detE$ is a supporting divisor of an extremal face, $F$; let $R_i$ the
extremal rays spanning this face. Fix one of this ray, say $R = R_i$ and
let $\pi:X\rightarrow W$ be the elementary contraction associated to $R$.

We have $l(R)\geq n-2$; this implies first that the
exceptional loci are
disjoint if $n > 3$, proposition
(\ref{n=4}). Secondly, by the inequality (\ref{diswis}),
 we have
$$dimE(R)+dimF(R)\geq 2n - 3.$$
Therefore $dimE(R)=n-1$ and either $dimF(R) = n-1$ or $dimF(R) = n-2$;
if $Z := \rho (E)$ and $D=E(R)$ this implies that either $dimZ = 0$ or $1$.

If $dimZ = 1$ then $dim F(\pi) = n-2$ for all fibers (note that since the
contraction $\pi$ is elementary there cannot be fiber of dimension
$(n-1)$); thus we can apply
proposition
(\ref{bd}) with $r = (n-2)$. This will give the case 3-(iii) of the theorem.

Consider again the construction in section (\ref{tech}),
 in particular we refer to
the diagram (\ref{dia1}). Let $S$ be the extremal ray contracted
by $\varphi$; note that $l(S)\geq n-2$
and that the inequality (\ref{diswis}) gives
$$dimE(S)+dimF(S)\geq 3n - 6;$$
in particular, since $dim F(S) \leq dim F(R) $,
we have two cases, namely $dimE(S) = 2n-5$ and $dimF(S) = (n-1)$ or
$dimE(S) = 2n-4$ and $dimF(S) = (n-1)$ or $(n-2)$.

The case in which $dimE(S) = 2n-5$ will not occur. In fact, after "slicing",
(see \ref{adj}),
we would obtain a map $\varphi^{\prime}=\varphi_{|Z}$ which would be a small
contraction supported by a divisor of the type $K_Z+(n-2)L$
but this is impossible by the classification of \cite[Th 4]{Fu1}
(see also \cite{An}).

\medskip
Hence $dimE(S)=2n-4$, that is also $\varphi$ is divisorial.

Suppose that the general fiber of $\varphi$, $F(S)$, has
dimension $(n-2)$. After slicing we obtain a map
${\varphi}^{\prime}=\varphi_{|Z}: Z \ra T$
supported by $K_Z+(n-2)L$, where $L={\taut}_{|Z}$.
This map contracts divisors $D$ in $Z$ to curves; by
(\cite[Th 4]{Fu1})
we know that every fiber $F$ of this map is $\P^{(n-2)}$ and that
$D_{|F} = \O(-1)$ (actually this map is a blow up of a smooth curve in a
smooth variety).
In particular there are curves in $Y$, call them
$C$, such that $-E(S).C = 1$. We will discuss this case in
a while.

Suppose then the general fiber of $\varphi$, $F(S)$, has
dimension $(n-1)$; therefore all fibers have dimension $(n-1)$.
Slicing we obtain a map
${\varphi}^{\prime}=\varphi_{|Z}: Z \ra T$
supported by $K_Z+(n-2)L$, where $L={\taut}_{|Z}$.
This map contracts divisors $D$ in $Z$ to points; by (\cite{Fu1})
we know that these divisors are either $\P^{(n-1)}$
with normal bundle $\O(-2)$ or $\Q^{(n-1)}\subset \Proj^n$ with normal bundle
$\O(-1)$.
In the latter case we have as above that there are curves $C$ in $Y$,
 such that $-E(S).C = 1$.

In these cases observe that $E(S)\cdot \tilde{C}=0$, where $\tilde{C}$
is a curve in the fiber of $p$. Hence $E(S)=p^*(-M)$ for some $M\in Div(X)$.
Let $l$ be an extremal curve of $E(S)$. Then, by projection formula, we have
$-1=E(S)\cdot l=-M\cdot mC$ and thus  $M$ generates $Im[Pic(X)\ra Pic(D)]$,
hence
$M$ is $\pi$-ample;
note that in general it does not generate $Pic(D)$.
We study now the Hilbert polynomial of $M_{|D}$ to show
that $\Delta(D,M_{|D})=0$, where $\Delta(X,L)$ is
Fujita delta genus.
Let $\O_D(-K_X) \simeq\O_D(pM)$,
where $p=l(R)\geq n-2$, and $\O_D(-D)\simeq\O_D(qM)$
for some $p,q\in {\bf N}$. By adjunction
formula $\omega_D\simeq \O_D(-(p+q)M)$. By \cite[Lemma 2.2]{Ando} or
\cite[pag 179]{BS}, Serre duality and relative vanishing
we obtain that $q\leq 2$, the Hilbert polynomial is
$$P(D,M_{|D})= \frac{a}{(n-1)!}(t+1)\cdots(t+(n-2))(t+c)$$
and the only possibilities are $a=1, c=n-1, q=1 or 2$
and $a=2, c=(n-1)/2, q=1$. In particular $\Delta(D,M_{|D})=0$ and,
by Fujita classification, $D$ is equal to $\P^{(n-1)}$
or to $\Q^{(n-1)}\subset \Proj^n$. Now the rest of the claim in
3) i) and ii) follows easily.

It remains the case in which
${\varphi}^{\prime}=\varphi_{|Z}: Z \ra T$
contracts divisors $D= \P^{(n-1)}$
with normal bundle $\O(-2)$ to points.
We can apply the above proposition (\ref{fujita})
and show that the singularities of $W$ are the same as those of
$T$. Then, as in (\cite{Mo1}), this means that we can
factorize $\pi$ with the blow up of the singular point.
Let $X^{\prime}=Bl_{w}(W)$, then we have a birational map $g:X\rightarrow
X^{\prime}$. Note that $X^{\prime}$ is smooth and that $g$ is finite.
Actually
it is an isomorphism outside $D$ and cannot contract any curve of
$D$. Assume to the contrary that $g$ contracts a curve $B\subset D$;
 let $N\in Pic(X^{\prime})$ be an ample divisor then we have
$g^*N\cdot B=0$ while $g^*N\cdot C\not=0$ contradiction.
Thus by Zarisky's main theorem $g$ is an
isomorphism. This gives a case in 3)i).

\end{document}